\documentclass[a4paper,british]{amsart}
\usepackage{babel}
\usepackage[T1]{fontenc}
\usepackage{amscd}
\usepackage{url}
\usepackage{amssymb}

\newtheorem{thm}{Theorem}[section]
\newtheorem{cor}[thm]{Corollary}
\newtheorem{prop}[thm]{Proposition}

\theoremstyle{remark}
\newtheorem*{rmk}{Remark}

\theoremstyle{definition}

\newcommand{\ie}{\textit{i.e. }}
\newcommand{\p}{\mathbb{P}}

\newcommand{\F}{\mathcal{F}}

\newcommand{\Gr}{\mathbb{G}}

\DeclareMathOperator{\length}{length}

\DeclareMathOperator{\red}{red}

\begin{document}

\title{Congruences of Lines with one-dimensional focal locus}

\author{Pietro De Poi}
\thanks{This research was partially supported by the DFG Forschungsschwerpunkt 
``Globalen Methoden in der Komplexen Geometrie'', and the EU (under the EAGER 
network).}
\address{Mathematisches Institut\\  
Universit\"at Bayreuth\\ 
Lehrstuhl VIII \\
Universit\"atsstra\ss e 30\\
D-95447 Bayreuth\\ 
Germany\\
}
\email{depoi@mathsun1.univ.trieste.it}

\date{\today}

\begin{abstract}
In this article we obtain the classification of the congruences of lines  
with one-dimensional focal locus. It turns out that one can restrict to 
study the case of $\p^3$. 
\end{abstract}

\maketitle

\bibliographystyle{amsalpha}

\section*{Introduction}
A congruence $B$ of lines in $\p^n$ is a family of lines of dimension $n-1$. The
order of $B$ is the number of lines of $B$ passing through a general point in $\p^n$ and 
its class is the number of lines of the congruence contained in a general 
hyperplane and meeting a general line 
contained in it.

The study of congruences of lines in $\p^3$ was 
started by E.~Kummer in \cite{K}, in which he gave a classification of those of 
order one. 
More recently, congruences of lines in 
$\p^3$ were studied in \cite{G} by N. Goldstein, who 
tried to classify these from the point of view of their focal locus. The focal locus is, 
roughly speaking, the set of critical values of the natural map between the total space $\Lambda$ 
of the family $B$ and $\p^n$. 
Successively, Z. Ran in \cite{R} studied the surfaces of order one in the 
general Grassmannian $\Gr(r,n)$ \ie families of $r$-planes in $\p^n$ 
for which the general $(n-r-2)$-plane meets only one element of the family. 
He gave a classification of such surfaces, in particular, in the case of 
$\p^3$, obtaining a modern and more correct proof of Kummer's 
classification of first order congruences of lines. Moreover, he proved that if the class of 
these congruences in $\p^3$ is greater 
than three, these are not smooth, as conjectured by I. Sols. Another proof 
of Kummer's classification is given by F. L. Zak and others in \cite{Z}. 
Moreover, E. Arrondo and I. Sols, in \cite{Ar} classified all the smooth 
congruences with small invariants. Successively, E. Arrondo and M. Gross in 
\cite{AG} classified all the smooth congruences in $\p^3$ with a 
fundamental curve and this classification was extended to all smooth 
congruences with a fundamental curve by E. Arrondo, M. Bertolini and 
C. Turrini in \cite{ABT}.

In this paper the classification of congruences with focal locus of dimension 
one (\ie if the focal locus reduces to be only a---possibly reducible or 
non-reduced---fundamental curve) is obtained. 
In Section~\ref{sec:1} we fix the notation and we recall some known results. 
In Section~\ref{sec:2} we obtain some general results on the possible minimal 
dimension of the focal locus of a congruence in $\p^n$: first of all, we prove in 
Theorem~\ref{thm:cata} that 
if the focal locus has codimension at least two, then the congruence has order 
either zero or one. This result has been suggested to us by F. Catanese. 
In Theorem~\ref{thm:dim} 
we prove that the dimension of the focal locus 
is at least $\frac{n-1}{2}$---if the congruence is 
not a star of lines, \ie the set of lines passing through a point $P\in \p^n$ (in which case 
the focal locus has support in $P$). Moreover, we give also some results concerning the case 
of dimension  $\frac{n-1}{2}$. 

In Section~\ref{sec:3} we prove the main result of this paper, namely Theorem~\ref{thm:int}; 
in order to state it, we need the following notation: first of all,  
let $\ell$  be a fixed line in $\p^3$ and $\p^1_\ell$ the set
of the planes containing $\ell$. Let $\phi$ be 
a general nonconstant morphism from $\p^1_\ell$ to $\ell$ and let $\Pi$ be a general
element in $\p^1_\ell$.
We define $\p^1_{\phi(\Pi),\Pi}$ as the pencil of lines passing through the point $\phi(\Pi)$ and
contained in $\Pi$.

\begin{thm}\label{thm:int}
If the focal locus $F$ of a congruence of lines $B$ in $\p^n$ 
has dimension one, then $n=2,3$. If $n=2$ $F$ is an irreducible curve of 
degree greater than one and $B$ is its dual; if $n=3$ $B$ is a first 
order congruence, and we have the following possibilities: 
\begin{enumerate}
\item $F$ is an irreducible curve, which is either: 
\begin{enumerate} 
\item\label{1st} a rational normal curve $C_3$ in $\p^3$. In this case $B$ 
is the family of secant lines of $C_3$ and has bidegree $(1,3)$ in the Grassmannian; 
\item\label{punto} or a non-reduced scheme whose support is a line $\ell$, 
and the congruence is---using the above
notation---$\cup_{\Pi\in \p^1_\ell}\p^1_{\phi(\Pi),\Pi}$; 
if we set $d:=\deg(\phi)$, 
the congruence has bidegree $(1,d)$; or 
\end{enumerate}
\item\label{2nd} $F$ is a reducible curve, 
union of a line $F_1$ and  a rational curve $F_2$ such that 
$\length(F_1\cap F_2)=\deg(F_2)-1$; 
$B$ is the family of lines meeting $F_1$ and $F_2$ and it  
has bidegree $(1,\deg(F_2))$. 
\end{enumerate} 
\end{thm}

\subsection*{acknowledgments}
This work has been quite improved by many valuable suggestions of Prof. F. Catanese.  

I would also like to thank Prof. E. Arrondo, Prof. E. Mezzetti for useful discussions on
these topics and Prof. F. Zak for some references, especially
\cite{Z}.


\section{Notation and Preliminaries}\label{sec:1}

We will work with schemes and varieties over the complex field $\mathbb C$, 
with the standard definitions and notation as in \cite{H}.
For us a \emph{variety} will always be projective. 
More information about general results and references 
about families of lines, focal diagrams and congruences can be found in 
\cite{DP2}, \cite{DP3} and \cite{DPi}, or \cite{PDP}. 
Besides, we refer to \cite{GH} for notation about 
Schubert cycles and to \cite{Fu} for the definitions and results of 
intersection theory. So we denote  by $\sigma_{a_0,a_1}$ 
the Schubert cycle of the lines in $\p^n$ contained in a fixed  
$(n-a_1)$-dimensional subspace $H\subset\p^n$ and which meet a fixed $(n-1-a_0)$-dimensional 
subspace $\Pi\subset H$. 
Here we recall that a \emph{congruence of lines} in $\p^n$ 
is a (flat) family $(\Lambda,B,p)$ of lines in $\p^n$ obtained 
as the pull-back of the universal family under the 
desingularization map of a subvariety $B'$ of dimension $n-1$ of the Grassmannian 
$\Gr(1,n)$ of lines in $\p^n$. So $\Lambda\subset B\times \p^n$ and  
$p$ is the restriction of the projection 
$p_1:B\times\p^n\rightarrow B$ to $\Lambda$, while we will denote the 
restriction of $p_2:B\times\p^n\rightarrow\p^n$ by $f$. 
$\Lambda_b:=p^{-1}(b)$, $(b\in B)$ will be a line of the family and 
$f(\Lambda_b)=:\Lambda(b)$ is the corresponding line in $\p^n$. 
A point $y\in \p^n$ is called 
\emph{fundamental} if its fibre $f^{-1}(y)$ has dimension greater than the 
dimension of the general one. The \emph{fundamental locus} is the set of the 
fundamental points. It is denoted by $\Psi$. Moreover the locus of the points  
$y\in \p^n$ for which the  fibre $f^{-1}(y)$ has positive dimension 
will be denoted by $\Phi$. 

Since $\Lambda$ is smooth of dimension $n$, we define the \emph{focal divisor} 
$R\subset\Lambda$ as the ramification divisor of $f$. The schematic image 
of the focal divisor $R$ under $f$ (see, for example, \cite{H}) is the 
\emph{focal locus} $F\subset \p^n$. 
Clearly, we have $\Psi\subset\Phi\subset (F)_{\red}$. 

To a congruence is associated a \emph{sequence of degrees} 
$(a_0,\dotsc,a_\nu)$ since $B$ is rationally equivalent  to  
a linear combination of the 
Schubert cycles of dimension $n-1$ in the Grassmannian $\Gr(1,n)$: 
$$
[B]=\sum_{i=0}^{\nu}a_i\sigma_{n-1-i,i}
$$ 
(where we set $\nu:=\left [\frac{n-1}{2}\right ]$); in particular, the \emph{order} $a_0$ 
is the number of lines of $B$ passing through a general point in $\p^n$, and 
the \emph{class} $a_1$ is the number of lines intersecting a general 
line $L$ and contained in a general hyperplane $H\supset L$ 
(\ie as a Schubert cycle, $B\cdot \sigma_{n-2,1}$). $a_0$ is the degree of $f:\Lambda\rightarrow\p^n$: 
thus if $a_0>0$, 
$\Psi=\Phi$, while if $a_0=1$, set-theoretically,  $\Psi=\Phi=F$. 

An important result---independent of order and class---is the following:
\begin{prop}[C. Segre, \cite{sg}]\label{prop:fofi}
On every line $\Lambda_b\subset \Lambda$ of the family, 
the focal divisor $R$ 
either coincides with the whole  $\Lambda_b$---in which case $\Lambda(b)$ 
is called a \emph{focal line}---or is a zero dimensional 
subscheme of $\Lambda_b$ of length $n-1$. Moreover, in the latter case, 
if $\Lambda$ is a first order congruence, $\Psi=F$ and 
$F\cap \Lambda(b)$ has length  $n-1$. 
\end{prop}

This result was proven classically in  \cite{sg}; the first modern proof in the case of the congruences 
(\ie families of dimension two) of planes in $\p^4$ is in \cite{CS}. See \cite{DP2}, Proposition~1 
for a proof of Proposition~\ref{prop:fofi}. 
 

\section{Congruences of lines in $\p^n$}
\label{sec:2}

We start with the following result, suggested to us by F. Catanese:

\begin{thm}\label{thm:cata}
Let $B$ be a congruence whose focal locus $F$ has codimension at least two. 
Then, $B$ has order either zero or one. 
\end{thm}

\begin{proof}
Let us consider the restriction of the map $f:\Lambda\rightarrow\p^n$ to the 
set $\Lambda\setminus f^{-1}(F)$. Then, either  $f^{-1}(F)=\Lambda$, in which 
case $B$ is a congruence of order zero, or the map $f\mid_{\Lambda\setminus f^{-1}(F)}$
defines an unramified covering of the set $\p^n\setminus F$. But it is a 
well-known fact that---by dimensional reasons---$\p^n\setminus F$ is 
simply connected and $\Lambda\setminus f^{-1}(F)$ is connected. Therefore, 
$f\mid_{\Lambda\setminus f^{-1}(F)}$ is a homeomorphism, hence $f$ is a birational 
map and $B$ is a first order congruence. 
\end{proof}

A central definition, introduced in \cite{DP2}, 
is the notion, $\forall d$ such that $0\le d\le n-2$, of \emph{fundamental $d$-locus}.  
This locus $\F^d$,  
is so defined: let $R_1,\dotsc,R_s$ be the \emph{horizontal components} 
(\ie $p_{\mid R_i}:R_i\rightarrow B$ is dominant).  
Now, let $F_i:=f(R_i)$, (as before, we take the schematic image); 
then, $\forall d$ such that $0\le d\le n-2$ we define $\F^d:=\cup_{\dim(F_i)=d}F_i$.  
By dimensional reasons, 
the lines of the congruence passing through $P_d\in \F^d$ form a family of dimension 
(at least) $n-1-d$ (so, set-theoretically, $\F^d\subset\Phi$). 

From this, we infer for example that if the focal locus has 
dimension zero, then the congruence is a star of lines (see Corollary~1 of 
\cite{DP2}). In this article, we will be interested in the subsequent case, 
\ie the case in which $\dim(F)=1$.  

We first note that if our congruence has order one, then there exists a $d$ such that 
$\F^d\neq \emptyset$.

The first important results of 
this paper are the following

\begin{thm}\label{thm:dim}
Let $B$ be a congruence of lines in $\p^n$. If $\dim(F):=i>0$,  
then $\frac{n-1}{2}\le i \le n-1$. 
Besides, if $i=\frac{n-1}{2}$, $(F)_{\red}$ 
is 
irreducible and the general line of the congruence meets 
$F$---set-theoretically---in only one point, 
then $F$ is---set-theoretically---an $i$-plane.
\end{thm}

\begin{proof} 
Clearly, we have only to prove that $\dim(F) \ge \frac{n-1}{2}$. 
Moreover, by Theorem \ref{thm:cata}, we 
can reduce to study the cases of order zero and one. 

If the congruence has order zero, $F=f(\Lambda)$ and $F$ 
contains a family of lines of dimension $n-1$. Then the  
the observation that  among  
varieties  of fixed dimension $i$, the projective space is the unique one 
which contains a family of lines of dimension greater than or equal to $2(i-1)$. 

Let now suppose that the order is one. Let $\F^i$ be the fundamental $i$-locus of 
maximal dimension $i(>0)$.  
Given $i+1$ general hyperplanes $H_0,\dotsc,H_i$ in $\p^n$ and the 
corresponding 
hyperplane sections of $\F^i$, $D_0,\dotsc,D_i$, the lines 
of $\Lambda$ which meet $D_j$ form a family of dimension $n-2$, $j=0,\dotsc,i$ 
which will fill a hypersurface $M_{D_j}$ in $\p^n$. Since 
$D_0\cap \dotsb\cap D_i=\emptyset$ and 
$M:=M_{D_0}\cap \dotsb \cap M_{D_i}\subset \p^n$ is such that $\dim(M)\ge n-1-i\ge 0$, take   
$Q\in M$. By definition, there exist $\ell_i\in D_i$ with $Q\in\ell_i$. If the $\ell_i$'s are 
different, then $Q\in F$. If all $\ell_i$'s are equal, it is absurd. It follows that $M\subset F$.


Now  $M\subset F$;  if  $\dim(F)<\frac{n-1}{2}$, 
then, since $\dim(M)\ge n-1-i$, we have 
$2i> n-1$, which cannot be, since $\F^i\subset F$.  

Let us now suppose that the equality holds, $D:=(F)_{\red}$ is 
irreducible, 
and the general line of the congruence meets 
$D$ in only one point.  
If $P,Q\in D$ are two general points, the set of lines of $\Lambda$ 
passing through them generate two cones $M_P,M_Q\subset\p^n$ 
of dimension $i+1$. 
Therefore, by our hypotheses $\dim(M_P\cap M_Q)\ge 1$ and $M_P\cap M_Q$ 
contains the line joining $P$ and $Q$. 
So, $D$ is a linear space. 
\end{proof}

\begin{prop}\label{prop:dim}
Let $\Lambda$ be a congruence in $\p^n$ and $D$  the component 
of the focal locus of maximal dimension $i>0$; if $\F^k$ is a fundamental 
$k$-locus contained in $D$, then $n-1-i\le k\le i$.
In particular, if $i=\frac{n-1}{2}$ we have only the $i$-fundamental 
locus $\F^i$.
\end{prop}

\begin{proof}
With notation as in the proof of the preceding theorem, 
we consider the hyperplane sections of $F^k$,  $D_0,\dotsc,D_k$, 
for which we have  $D_0\cap \dotsb\cap D_k=\emptyset$, 
and so $M_{D_0}\cap \dotsb \cap M_{D_k}\subset F$.
Therefore our thesis easily follows.
\end{proof}

\begin{cor}\label{cor:1}
Let $\Lambda$ be a congruence in $\p^n$; if $\dim(F)=1$, then $n\le 3$.  
\end{cor}

\begin{rmk}
The preceding results are new in the 
literature, although the method of proof of  Theorem~\ref{thm:dim} and Proposition~\ref{prop:dim}
is very similar to the respective proofs of Theorem~6 and Corollary~2 of \cite{DP2}. Actually, 
in Theorem~6 of \cite{DP2} we forgot the hypothesis---in the border case 
$i=\frac{n-1}{2}$---that the general line should meet the focal locus only 
once. 
\end{rmk}


\section{Congruences in $\p^3$}
\label{sec:3}

Our aim in this article is to classify the congruences with $\dim(F)=1$. 
In view of Corollary~\ref{cor:1}, it is sufficient to consider congruences 
in $\p^2$ and $\p^3$. The case of $\p^2$ is well known: all the congruences 
but the ones of order one (\ie the pencils) have a focal curve, and the 
congruence is given by the tangents of the focal curve. 
Therefore, we will study congruences in $\p^3$ with $\dim(F)=1$. By Theorem~\ref{thm:cata}, 
these are the first order congruences (but the star of lines) of $\p^3$, since 
the order cannot be zero; this can be deduced easily from the fact that 
the only surface containing a family of dimension two of lines is the plane. 
Actually, the first order congruences of $\p^3$ are classified in \cite{R} and 
\cite{Z}. Here we give another easy and---as far as we know---new proof of this classification. 
 
In what follows, we will also call the 
focal locus \emph{fundamental curve}.


\subsection{Congruences with Irreducible Fundamental Curve}

\begin{prop}\label{prop:irr}
If the fundamental curve $F$ of a congruence of lines in $\p^3$ 
is irreducible, then we have the following possibilities: 
\begin{enumerate}
\item if $F$ is reduced, we are in case~\eqref{1st} of 
Theorem~\textup{\ref{thm:int}};
\item  If $F$ is non-reduced, we are in case~\eqref{punto} of 
Theorem~\textup{\ref{thm:int}}.
\end{enumerate}
\end{prop} 

\begin{proof}
First of all, we consider the case of the secant lines of a curve, and 
we know that we have a first order congruence. 
We could apply Theorem~2.16 of \cite{DPi} to obtain that $\deg(F)=3$, but 
for the sake of completeness, we will prove this result here.  
Clearly we have that $\deg F\ge 3$ for degree reasons. Then, let us denote by 
$a$ the class  of our congruence, 
and by $r$ a general line in $\p^3$. Let us also denote by $V_r$ the scroll of
the lines of $B$ meeting $r$. An easy application of the Schubert calculus
shows that $\deg(V_r)=a+1$.

Let $V_r$ and $V_{r'}$ be two such surfaces; we denote by $k$ the
(algebraic) multiplicity of $F$ in $V_r$ (and $V_{r'}$).
First of all, we note that
the complete intersection of two general surfaces  $V_r$ and  $V_{r'}$
is a (reducible) curve $\Gamma$ whose
components are the focal locus $F$ and $(a+1)$ lines, \ie the lines
of $B$ meeting $r$ and $r'$: in fact, a point of $\Gamma$ not belonging
to the lines meeting  both $r$ and $r'$ is a focal point, since through it
pass infinitely many lines of $B$, \ie a line meeting $r$ and another one
meeting $r'$, and $r$ and $r'$ are general. 
The lines of $B$ meeting $r$ and $r'$
are $a+1$ since this is the degree of $V_r$. 

Then, by B\'ezout, we have that
\begin{align*}
\deg(V_r\cap V_{r'}) &=(a+1)^2\\
 &= a+1 + k^2d
\end{align*}
---where $d:=\deg(F)$---and we obtain
\begin{equation}\label{eq:cu1}
k^2d=(a+1)a.
\end{equation}

Besides, since a line of the congruence not belonging to the $(a+1)$ lines
meeting $r$ and $r'$  must intersect the curve $F$ in (at least) two points,
we deduce
\begin{equation}\label{eq:cu2}
a+1=2k.
\end{equation}
From formulas \eqref{eq:cu1} and \eqref{eq:cu2} we conclude
\begin{equation*}
k^2(4-d)-2k=0,
\end{equation*}
from which we get the result, \ie $d=3$.

So the only possibility is to have the 
twisted cubic $C^3$, 
which has, in fact, an apparent double point.
Since this curve has degree three, the bidegree is $(1,3)$.

Then we pass at the other case: by Theorem~\ref{thm:dim}, $(F)_{\red}=:\ell$ is a line;  
so, if we restrict the congruence to a (general) plane $\Pi$ through $\ell$, 
we get a congruence of lines in $\p^2$ such that its focal locus is contained 
in $\ell$, and therefore 
it is a pencil of lines with its fundamental point $P_\Pi$ in $\ell$. 
Then $\varphi$ of Theorem~\textup{\ref{thm:int}}, \eqref{punto} 
is the map defined by $\varphi(\Pi)=P_\Pi$.
This congruence has bidegree $(1,d)$, as it is easily seen if we consider
the lines of the congruence contained in a general plane.
\end{proof}


\subsection{Congruences with Reducible Fundamental Curve}

\begin{prop}\label{prop:red}
If the fundamental curve $F$ of a first order congruence of lines in $\p^3$ is 
reducible, then the congruence is as in Theorem~\textup{\ref{thm:int}}, 
case~\eqref{2nd}. 
\end{prop}

\begin{proof}
Let us denote by $Z$ the 0-dimensional scheme given by 
$F_1\cap F_2$ and we set $u=\length(Z)$.

Let $P$ be a general point in $\p^3$; then, if we set $\deg(F_1):=m_1\ge m_2=:\deg(F_2)$,  
the cone given by the join
$PF_j$ has degree $m_j$; as usual, by B\'ezout
$$
\deg(PF_1\cap PF_2)=m_1m_2.
$$
Since we have a congruence of order one, we obtain:
\begin{equation}\label{eq:m-1}
u=m_1m_2-1; 
\end{equation}
this is due to the fact that only one of the lines of $PF_1\cap PF_2$ is a line of the congruence; therefore 
the others must be the $u$ lines  of the join $PZ$. 

If $Q$ is a general point of $F_1$, there will pass 
$m_2(m_1-1)$ secant lines of $F_1$ through $Q$ meeting $F_2$ also; 
but these lines will pass through the points of $Z$, 
since if one of these lines did not intersect $Z$, this line would be a
focal line, and varying these lines when we vary $Q$ on $F_1$, 
we would obtain a focal surface. Besides, since $Q$ is general, we can 
suppose that it does not belong to the tangent cones of the two curves 
at the points of $Z$. Then we have that 
\begin{equation}\label{eq:mm}
u=(m_1-1)m_2,
\end{equation}
and by equation \eqref{eq:m-1} we obtain $m_2=1$, $u=m_1-1$. To see that $F_1$ is a  
rational curve, we can simply project it from a general point of the line $F_2$
onto a plane: by the Clebsh formula, the projected 
curve has geometric genus zero.

The class is---as usual---calculated intersecting with a general
plane.
\end{proof}



\subsection{Final Remarks about First Order Congruences} 
We finish this article analysing which of these congruences are smooth 
as surfaces in 
the Grassmannian $\Gr(1,3)$. 

\begin{prop}\label{thm:last}
The smooth first order congruences of lines $B$ in $\p^3$ are
\begin{enumerate}
\item the secants of the twisted cubic, in which case $B$ is the 
Veronese surface;
\item the join of a line and a smooth conic meeting in a point, in 
which case $B$ is a rational normal cubic scroll;
\item the join of two---non-meeting---lines in 
which case $B$ is a quadric.
\end{enumerate}
\end{prop}

\begin{proof}
We see from \cite{AG} and \cite{ABT} that the only possible bidegrees 
for a smooth congruence in $\p^3$ with a fundamental curve are $(1,1)$, $(1,2)$, $(1,3)$ 
and $(3,6)$. The last case can be excluded because has order three, by Theorem~\ref{thm:int}. 
From the description of these congruence given in \cite{Ar} we get 
the proposition. 
\end{proof}

\begin{rmk}
Observe that, if the fundamental curve is non-reduced, its support is a 
line $\ell$. Its image in $\Gr(1,3)$ is a singular point $L$ for $B$: more 
precisely, $B$ is a $2$-dimensional linear section,  tangent to $\Gr(1,3)$ 
at $L$. 

We observe finally that the case of two lines is the one in which $B$
is the intersection of the two tangent hyperplane sections of the 
Grassmannian corresponding to the two lines, 
\ie $B$ is a general $2$-dimensional linear section of $\Gr(1,3)$. 
From this description we see that the case in which the fundamental curve is a
line only is a limit case of this one, when the two lines coincide.  
\end{rmk}


\providecommand{\bysame}{\leavevmode\hbox to3em{\hrulefill}\thinspace}
\providecommand{\MR}{\relax\ifhmode\unskip\space\fi MR }
\providecommand{\MRhref}[2]{%
  \href{http://www.ams.org/mathscinet-getitem?mr=#1}{#2}
}
\providecommand{\href}[2]{#2}

\end{document}